\documentclass[a4paper]{article}
\usepackage[english]{babel}
\usepackage{a4wide}
\usepackage{amsmath,amssymb,amsthm}
\usepackage{mathtools}
\usepackage{ifpdf}
\usepackage[numbers,sort&compress]{natbib}

\ifpdf%
\usepackage[unicode,bookmarks=false,
hyperindex]{hyperref}
\hypersetup{pdfstartview=FitH,
linkbordercolor=[rgb]{0.5, 1, 1},
citebordercolor=[rgb]{0.5, 1, 0.5},
hypertexnames=false}
\else
\fi

\newcommand{\blfootnote}[2]{{
\renewcommand{\thefootnote}{#1}
\footnotetext[0]{#2}}}

\newtheorem{Theorem}{Theorem}
\newtheorem{Lemma}{Lemma}
\theoremstyle{remark}
\newtheorem{Remark}{Remark}

\title{Branching Random Walks in a Random Killing Environment with a Single Reproduction Source%
\blfootnote{$1$}{Department of Probability Theory, Lomonosov Moscow State University, Moscow, Russia}%
\blfootnote{$2$}{Department of Probability Theory and Mathematical Statistics, Steklov Mathematical Institute, Moscow, Russia}%
\blfootnote{$3$}{Laboratory of Stochastic Analysis and its Applications, National Research University Higher School of Economics, Moscow, Russia}%
\blfootnote{$4$}{Department of Mathematics and Statistics, University of North Carolina at Charlotte, Charlotte, NC, USA}%
}

\author{Vladimir Kutsenko$^{1,2}$, Stanislav Molchanov$^{3,4}$, and Elena Yarovaya$^{1,2}$}

\begin{document}
\date{}
\maketitle

\begin{abstract}We consider a continuous-time branching random walk on $\mathbb{Z}$  in a random non homogeneous environment. Particles can walk on the lattice points or disappear with random intensities. The process starts with one particle at initial time $t=0$. It can walk on the lattice points or  disappear with a random intensity until it reach the point, where initial particle
can split into two offspring. This lattice point we call reproduction source. The offspring of the initial particle evolve according to the same law, independently of each other and the entire prehistory. The aim of the paper is to study the conditions for the presence of exponential growth of the average number of particle at an every lattice point. For this purpose we investigate the spectrum of the random evolution operator of the average particle numbers. We derive the condition under which there is exponential growth with probability one. We also study the process under the violation of this condition and present the lower and upper estimates for the probability of exponential growth.

\medskip

\end{abstract}

\section{Introduction}

Apparently, stochastic processes which can be reduced to a model of branching random walks (BRWs) in a random environment first appeared in the field of statistical physics in a series of papers by Y.~Zeldovich and co-authors \cite{shukurov84,zeldovich85,zeldovich87}. These processes were further developed in the work of S.~Molchanov and J.~G\"{a}rtner, where the authors introduced basic concepts for  BRWs in a random environment and developed approaches to BRW analysis~\cite{GM90}.  The model studied in this work coincides with the previously introduced parabolic Anderson  model \cite{anderson}, and the paper~\cite{GM90} itself is recognized as fundamental and has sparked active research on applications of the Anderson model in various fields~\cite{konig16}.  

Subsequent papers \cite{M94,GM98, albeverio2000} considered various BRW models and the effects arising in them. In particular, the irregularity of the growth rate of particle number moments (intermittency) was investigated, and a method to study the spectrum of random operators occurring in BRW in a random environments was described. The intermittency was analyzed by studying the asymptotics of the total particle number moments averaged over the environment under the assumption of an asymptotically Weibull distribution of the right tail of the random potential, i.e. the difference between the splitting rate and the death rate. The Feyman-Katz representations for moments developed in the \cite{GM90, M94,GM98, albeverio2000} allowed to obtain a number of results concerning the asymptotic description of the averaged characteristics of the BRWs, see, e.g., \cite{yarovaya2012, gun13, gun15, kutsenko23}. 

In addition to particle number moments averaged over the environment, other characteristics were also considered in papers about BRWs in a random environment. For example, E.~Chernousova with co-authors \cite{chernousova22} investigated the probability of BRW survival in a random environments and proved the absence of a steady state. The present work also poses a relatively new problem in the field of BRW in non homogeneous  random environments. The purpose of this paper is to investigate the probability of exponential growth of mean particle numbers for the BRW model on a one-dimensional lattice $\mathbb{Z}$ in a random killing environment with a single reproduction source. Note that a part of the results for this model is announced in \cite{umn}.

\section{Model Description}

Let us consider a branching random walk (BRW) on a one-dimensional 
lattice~$\mathbb{Z}$ with continuous time. On the lattice we define a field of 
of independent identically distributed random variables $\mathcal M = \left\{ 
\mu(x,\cdot), x\in \mathbb{Z} \backslash \{0\} \right\}$, which are defined 
on some probability space $\left(\Omega,\mathcal {F}, 
\mathbb{P}\right)$. We assume that each random variable 
$\mu(x,\cdot)$ takes values from the closed interval $[0,c]$ with $c\geqslant0$ and is a mixture of discrete and absolutely continuous random variables (r.v.). Also we assume that the continuous component of $\mu(x,\cdot)$ has a positive density on $[0,c]$. The field $\mathcal M$ forms on $\mathbb{Z}$
a ``random killing environment'' that determines the intensity 
of particle death in the BRW. In addition, we introduce the parameter $\Lambda 
\geqslant 0$, which is responsible for the intensity of particle multiplication at zero, and 
parameter $\varkappa > 0$, which controls the intensity of particle walking on the lattice.

Suppose that at time $t=0$ there is a realization of the $\mathcal M$ field denoted by 
$\mathcal M (\omega) = \left\{ \mu(x,\omega), x\in \mathbb{Z}\backslash 
\{0\}, \omega \in \Omega \right\}$. 
Also assume that the process at time $t=0$ starts with a single particle at some point $x \in \mathbb{Z}$.
The further evolution proceeds as follows. If the particle is at 
zero, then in time $h \rightarrow 0$ with probability $\Lambda h + o(h)$  it splits in two, with probability $\varkappa h + o(h)$ it moves equally likely
to one of the neighboring points, and with remaining probability $1 - \Lambda h - 
\varkappa h + o(h)$ it remains in place.  If the particle is at a point 
$ x \neq 0$, then in time $h \rightarrow 0$, with probability $\mu(x,\omega) h 
+ o(h)$ it disappears, with probability $\varkappa h + o(h)$ it moves
to one of the neighboring points, and with remaining probability $1 - \mu(x,\omega)  
- \varkappa h + o(h)$ it stays in place. The new particles evolve according to the same law independently of each other and of all prehistory.

The introduced process is Markovian, and can be described in terms of a set of exponential and polynomial variables. This description may be more convenient for the perception of the model.
Let us introduce the average waiting time $\tau(x)$ at an arbitrary point $x \in \mathbb{Z}$: 
$$\tau(x) = 
\begin{cases}
    \left(\varkappa + \Lambda\right)^{-1},& \text{if $x = 0$;} \\
    \left(\varkappa + \mu(x,\omega)\right)^{-1},& \text{if $x \not= 0$.}
\end{cases} $$
The evolution of a particle located at point $x$ is as follows. If the particle is at zero, it waits for an exponentially distributed time with parameter $\tau(0) ^{-1}$, and then instantly splits in two or moves equiprobably to one of the neighboring lattice points. The choice between these two events is made with corresponding probabilities $\left. \Lambda \tau(0)\right.$ and $\left. \varkappa\tau(0)\right.$. If the particle is at a point $x$ outside zero, it also waits an exponentially distributed time with parameter $\tau(x) ^{-1}$ and then vanishes instantaneously or moves equiprobably to one of the neighboring lattice points. The choice between these two events is made with corresponding probabilities $\left. \mu(x,\omega) \tau(x)\right.$ and $\left. \varkappa\tau(x)\right.$. The evolution of particles occurs independently of each other and of all prehistory.

We will show later on that the branching process 
of particles at the point $x \in \mathbb{Z}$ can be conveniently described by the potential $V(x,\omega)$ which reflects the criticality of the branching process at each point:
$$V(x,\omega) = \begin{cases}
\Lambda,& x = 0; \\
-\mu(x,\omega),& x \not= 0
\end{cases}
$$
or
\[
    V(x,\omega) = \Lambda \delta_0(x) -\mu(x,\omega) (1- \delta_0(x) ).
\]

The BRW at time $t$ due to Markov property  can be completely described by the set of particle numbers at time $t$ at the points $y\in \mathbb{Z}$ denoted by $N_{t}(y,\omega)$ . However, $N_{t}(y,\omega)$ is a random variable and hence is difficult to investigate. Therefore, it is common to consider the average 
particle number \cite{GM90,albeverio2000}: $$m_1(t,x,y,\omega) = \mathbb{E}_{x}  
N_{t}\left(y, \omega\right),$$ where $\mathbb{E}_{x}$ is the mathematical 
expectation under the condition that at time $t=0$ there is one particle at point $x$.

By $F_{\mu}$ we further denote an distribution function of $\mu(x)$. In this paper, we are interested in
the probability $P(\Lambda,\varkappa,F_{\mu})$ of the realization of an environments in which there is 
an exponential growth of $m_1(t,x,y,\omega)$ 
 for given parameters $\Lambda$, $\varkappa$, and $F_{\mu}$.  We will refer to such exponential growth as ``supercriticality''.
 Formal definition is as follows: 
 $$P(\Lambda,\varkappa,F_{\mu}) = \mathbb{P}\left\{ \omega \in \Omega : \exists \lambda,C(x,y)>0: \lim\limits_{t \rightarrow \infty} \frac{m_1(t,x,y,\omega) }{C(x,y)e^{\lambda t}}  = 1, \  \forall x,y \in \mathbb{Z} \right\},$$
 where $C(x,y) = C(x,y,\omega,\Lambda, \varkappa, F_{\mu}) $, $\lambda = \lambda(\omega,\Lambda, \varkappa, F_{\mu}) $.  Note that we require exponential growth of the average particle population simultaneously in all points of the lattice. However, further we will show that this condition is equivalent to exponential growth at least in one point. Intuitively speaking,  the exponential growth at one point ``is spreaded'' over the whole lattice with help of  random walk.

The purpose of this paper is to estimate $P(\Lambda,\varkappa,F_{\mu})$ as a function of $\Lambda$, 
$\varkappa$ and $F_{\mu}$. To do this, we first use the standard approach described, e.g., in~\cite{albeverio2000, yarovaya2012}, 
and write the Cauchy problem for $m_1(t,x,y, \omega)$:
\begin{equation}
\label{eq1}
\begin{aligned}
\frac{\partial m_1(t,x, y, \omega)}{\partial t} &= \left( \varkappa \Delta m_1 \right)(t,x, y, \omega)+V(x, \omega) m_1(t,x, y, \omega), \\
m_1(0, x, y) &= \delta_y(x),
\end{aligned}
\end{equation}
where $\varkappa\Delta f(x) = \frac{\varkappa}{2} \sum_{|x'- x| = 1} \left( f(x')-f(x) \right)$ is the discrete Laplace operator on $\mathbb Z$, and the sign $|\cdot|$ denotes the lattice distance on the $l_1$ norm. Here and below we assume that all 
operators are defined on $l_2(\mathbb{Z})$.

 Let us introduce a random self-adjoint 
operator $ H(\omega) = \varkappa \Delta + V(x, \omega)$ to rewrite the problem \eqref{eq1} in a simpler form
\begin{equation}
\label{eq1.1}
\begin{aligned}
\frac{\partial m_1(t,x, y, \omega)}{\partial t} &= H(\omega) m_1(t,x, y, \omega), \\
m_1(0, x, y) &= \delta_y(x).
\end{aligned}
\end{equation}
In problems of this kind, the behavior of $m_1$, and hence the quantity $P(\Lambda,\varkappa,F_{\mu})$, depends on the spectrum structure of the random operator $H(\omega)$. Therefore the present work is mainly devoted to the study of the spectrum of the  $H(\omega)$. In Sections~\ref{sec3} and~\ref{sec4}, it is shown that the spectrum of $\sigma(H(\omega))$ consists of a non-positive non-random part and can contain a positive random eigenvalue; in Section~\ref{sec5} we derive the condition under which $P(\Lambda,\varkappa,F_{\mu}) = 1$; we address violation of this condition in Sections~\ref{sec6} and~\ref{sec7}, where we present the lower and upper estimates for $P(\Lambda,\varkappa,F_{\mu})$. The main proofs are given in the text of the paper after the corresponding statements, while the auxiliary proofs are placed in Section~\ref{sec9}.

\section{The non-random part of the spectrum of the evolutionary operator} \label{sec3}

 We obtained  the results of this and the next section using the technique described in \cite{M94}. In these sections, we prove the results for the Cauchy problem in arbitrary dimension $d \in \mathbb{N}$, although the case $d = 1$ is sufficient to study our model. Consider the following Cauchy problem for $m_1(t,x,y, \omega)$:
\begin{equation}
\begin{aligned}
\frac{\partial m_1(t,x, y, \omega)}{\partial t} &= \left( \varkappa \Delta m_1 \right)(t,x, y, \omega)+V(x, \omega) m_1(t,x, y, \omega), \\
m_1(0, x, y) &= \delta_y(x),
\end{aligned}
\end{equation}
where $\varkappa\Delta f(x) = \frac{\varkappa}{2d} \sum_{|x'- x| = 1} \left( f(x')-f(x) \right)$ is the discrete Laplace operator on $\mathbb Z^d$, and the sign $|\cdot|$ denotes the lattice distance on the $l_1$ norm. 

For convenience of reasoning, let us introduce the averaging operator:
\[
    \left(\varkappa \bar\Delta f\right)(x) = \frac{\varkappa}{2d} \sum_{|x'- x| = 1} f(x'),
\]
where $\varkappa \bar\Delta f(x) = \frac{\varkappa}{2d} \sum_{|x'- x| = 1} f(x') $. 
The Laplace operator $\varkappa \Delta$ can be represented as the difference of the averaging operator and the multiplication operator:
\[
    \left(\varkappa \Delta f\right)(x) = \varkappa \bar\Delta f(x) - \varkappa f(x).
\]

Consider a simplified operator $H_\mu(\omega)$ for which the splitting intensity at zero is absent and the death intensity at zero $\mu(0,\omega)$ is defined in the same way as $\mu(x,\omega)$ for $x \in \mathbb{Z} \backslash \{0\}$. 
$$ H_\mu(\omega) = \varkappa \Delta - \mu(x, \omega) = \varkappa \bar\Delta - \varkappa - \mu(x, \omega).$$
The operator $H(\omega)$ can be viewed as a random one-point perturbation of the operator $H_\mu(\omega)$ at zero. Therefore, the essential spectra of these operators coincide \cite{akh_glazman}: $$\sigma_{ess}(H(\omega)) = \sigma_{ess}(H_\mu(\omega)).$$ Furthermore, a single-point perturbation can only produce at most one positive eigenvalue. Thus, the first problem is to investigate the essential spectrum of the operator $ H_\mu(\omega)$.

For convenience, we give the formulations of the lemmas from the works of \cite{akh_glazman, yarovaya2007}, which will be needed to study the spectrum of the operator $ H_\mu(\omega)$. 
\begin{Lemma}[see, e.g., \cite{akh_glazman}] \label{lemma1}
The number $\lambda$ belongs to the essential spectrum of the operator $H_\mu$ if we can construct a sequence of ``almost eigenfunctions'', i.e.:
\begin{equation}
  \begin{aligned}
 \exists\left\{f_n \in l_2 (\mathbb{Z}^d ): \quad\left\|f_n\right\|=1, \quad\left(f_n, f_m\right)=\delta(n,m),\right. 
  \left.\left\|H_\mu f_n- \lambda f_n\right\| \rightarrow 0, \quad n \rightarrow \infty\right\} .
\end{aligned}  
\end{equation}
\end{Lemma}

\begin{Lemma}[see, e.g., \cite{yarovaya2007}] \label{lemma2}
The spectrum of the operator $\varkappa \Delta$ is equal to $[-2\varkappa;0]$. For an eigenvalue $\lambda \in [-2\varkappa;0]$, there exists a representation 
\begin{equation*}
\lambda = \frac \varkappa {d} \sum_{i = 1}^d \cos(\phi_i) - \varkappa,
\end{equation*}
for some $\overrightarrow{\phi} = (\phi_1, \dots, \phi_d)$, $\phi_i \in [-\pi,\pi]$. The corresponding function
$ \psi_\lambda (x) = \exp \{ i(\overrightarrow{\phi}, x)\}$
is an eigenfunction for $\lambda$. As a consequence, the spectrum of the operator $\varkappa \bar\Delta$ is equal to $[-\varkappa;\varkappa]$.
\end{Lemma}
Using these lemmas and the proof scheme from \cite{M94}, we get the following result.

\begin{Lemma}\label{lemma3}
The spectrum of the operator $H_\mu(\omega)$ almost sure consists of only the essential part, which is equal to the interval $[-2\varkappa-c;0]$.
\end{Lemma}
\textit{Proof}.
The operator $H_\mu(\omega)$ is the sum of the averaging operator $\varkappa \bar \Delta$ and the multiplication by the function $-\mu(x,\omega)-\varkappa$.
Due to Lemma~\ref{lemma2}, the operator $\varkappa \bar \Delta$ has a spectrum equal to $[-\varkappa;\varkappa]$ and a norm equal to $\varkappa$. In turn, the spectrum of the multiplication operator on the function $-\mu(x,\omega)$ is equal to the closure of the set of values of this function. For almost sure (a.s.) any $\omega$, this closure is equal to the interval $[-c;0]$ by virtue of the definition of $\mu(x,\omega)$. Therefore, the spectrum of the combined operator $-\mu(x,\omega)-\varkappa$ is equal to $[-\varkappa-c;-\varkappa]$.

The operator $H_\mu(\omega)$ can be considered as a perturbation of the self-adjoint operator $-\mu(x,\omega)-\varkappa$ by the self-adjoint operator $\varkappa \bar \Delta$. In such a case, according to perturbation theory \cite{kato}, the spectrum of the operator $H_\mu$ will differ from the interval $[-\varkappa-c,-\varkappa]$ by at most $\varkappa$, leading to the following inclusion:
\begin{equation} \label{biba_1}
    \sigma (H_\mu) \subseteq [-2\varkappa-c;0].
\end{equation}

To show the inverse inclusion, we use Lemma~\ref{lemma1}, and for each $\lambda \in [-2\varkappa-c;0]$ we will construct a sequence of ``almost eigenfunctions'' $\{ f_n(x) \},$ $f_i(x) \in l_2 (\mathbb Z ^d)$. Note that we construct a sequence function for each fixed $\omega$, i.e., $\{ f_n(x) \}$ = $\{ f_n(x,\omega) \}$. 

Let us represent $\lambda$ as $\lambda = a + b$, $a \in [-2\varkappa;0]$ $b \in [-c;0]$. We need to construct $f_n$ such that the approximations $\varkappa \Delta f_n \approx a f_n$ and $-\mu(x,\omega) f_n \approx b f_n$ are true in some sense. Because then 
$$\left(\varkappa \Delta -\mu(x,\omega) \right) f_n \approx (a + b) f_n = \lambda f_n.$$

The condition $\varkappa \Delta f_n \approx a f_n$ requires that the function be ``almost everywhere'' similar to $\exp \{ i(\overrightarrow{\phi}, x)\}$ of Lemma~\ref{lemma2} with a suitable $\overrightarrow{\phi}$. The condition $-\mu(x,\omega) f_n \approx b f_n$ requires that the function be nonzero only on the region where $-\mu(x,\omega) \approx b$.

It turns out that the functions satisfying both conditions are indicators of the balls on which $-\mu(x,\omega) \in [b-\varepsilon; b+\varepsilon]$ for sufficiently small $\varepsilon>0$. Therefore $-\mu(x,\omega)\approx b$ and the multiplication operator will act ``almost eigen-like''. The mixing operator $\varkappa \Delta$ will act ``almost eigen-like'' inside and outside such balls, but not at the boundary. Therefore, the radius of the balls must increase to infinity so that the ``non-eigen'' action of $\varkappa \Delta$ tends to zero. The exact construction of the system of functions $\{ f_n \}$ and the proof they are ``almost eigenfunctions'' is given in Subsection~\ref{sec9_1}. 

In summary, for any $\lambda \in [-2\varkappa-c;0]$ we can construct a sequence of almost eigenfunctions functions $\{ f_n \}$, and hence
\begin{equation}  \label{boba_1}
\sigma(H_\mu)\supseteq [-2\varkappa-c;0].
\end{equation}
Inclusions \eqref{biba_1} and \eqref{boba_1} complete the proof of the lemma. \qed

Thus the result of the section is as follows: $$\sigma (H_\mu (\omega)) = \sigma_{ess} (H_\mu (\omega)) = [-2\varkappa-c;0].$$

\section{The random part of the spectrum of the evolutionary operator} \label{section_4} \label{sec4}

 Let us return to the operator $H(\omega) = \varkappa \Delta + V(x, \omega).$ As we have already mentioned, it can be viewed as a random one-point perturbation of the previously described operator $H_\mu$ with an essential spectrum $\sigma (H_\mu) = [-2\varkappa-c;0]$. By the Weyl criterion \cite{akh_glazman}, under compact perturbation the essential spectrum of the operator does not change, while one positive eigenvalue may appear, which we will denote by $\lambda(\omega)$: $$\sigma (H(\omega)) = [-2\varkappa-c;0] \cup \lambda(\omega).$$
 
As we mentioned under equation \eqref{eq1.1} the structure of the $\sigma (H(\omega))$ defines the behaviour BRW. In particular if $\lambda(\omega) >0$, an exponential growth of the average particle number is observed, see, e.g. \cite{yarovaya2007}. Thus, the study of the probability of exponential growth is reduced to the study of the probability of the appearance of a positive eigenvalue: 
$$P(\Lambda,\varkappa,F_{\mu}) = \mathbb{P} \left\{\omega: \exists \  \lambda(\omega) \in \sigma\left(H(\omega)\right): \lambda(\omega)>0 \right\}.$$

Let us formulate the problem of finding the eigenvalue of $\lambda(\omega)$ with the corresponding eigenfunction $u(x)$. Note that from   $u(0) = 0$ it follows that $u(x) \equiv 0$. Therefore, without restricting generality, let $u(0) = 1$: 
\begin{equation}
\label{eq3.1}
\begin{aligned}
& \left( \varkappa \Delta +V(x, \omega) \right)u(x) = \lambda u(x),\\
& u(0) =1.
\end{aligned}
\end{equation}
For convenience, equation \eqref{eq3.1} can be decomposed into two equations. When $x = 0$ it takes the following form:
\begin{equation}
\label{eq3.2}
\begin{aligned}
& \left( \varkappa \Delta + \Lambda - \lambda \right)u(0) = 0, \\
& u(0) =1.
\end{aligned}
\end{equation}
When $x \not= 0$, it takes the following form:
\begin{equation}
\label{eq3.3}
\begin{aligned}
& \left( \varkappa \Delta - \mu(x, \omega) - \lambda \right)u(x) = 0, \\
& u(0) =1.
\end{aligned}
\end{equation}
For simplicity of the formulas, we introduce the following notations:
$$ E = \varkappa +\lambda. $$
Due to this notation, equation \eqref{eq3.3} takes the following form:
\begin{equation}
\label{eq3.4}
\begin{aligned}
& \left( \varkappa \bar \Delta - ( \mu(x, \omega) + E) \right)u(x) = 0, \quad x \not=0,\\
& u(0) =1.
\end{aligned}
\end{equation}
Let us move on to finding the solution to this system of equations.
\begin{Lemma} \label{lemma_4}
    The solution to equation \eqref{eq3.4} when $x \not = 0$ is given by the following formula:
    \begin{equation}
    \label{u_formula}
        u(x) = \sum\limits_{\gamma: x \leadsto 0 } \prod\limits_{z \in \gamma} \left( \frac{\varkappa / 2d}{\mu(z,\omega) + E} \right),
    \end{equation}
    where by $\gamma: a \leadsto b = \{ a = x_1 ,\dots x_{n} \not= b\}$ we denote the path from point $a$ to point $b$ through the neighboring points of the lattice, and: a) the path does not intersect $0$; b) point $b$ is not considered to be included in the path $\gamma$.    
    The solution given by formula \eqref{u_formula} makes sense for any $\lambda >0$ in any dimension $d \in \mathbb{N}$.
\end{Lemma}
The first part of the lemma is verified by directly substituting formula \eqref{u_formula} into the problem \eqref{eq3.3}. The correct definiteness of the expression \eqref{u_formula} for any $\lambda > 0$ in dimension is checked by combinatorial reasoning and asymptotic methods. The full proof of Lemma~\ref{lemma_4} is given in Subsection~\ref{sec9_2}.

\begin{Remark}
Lemma~\ref{lemma_4} is a special case of the popular (especially in the physics literature) path expansion of the resolvent, but it is usually applied to the $\lambda$ from essential spectrum of Laplacian, i.e. $\lambda < 0$ in our case. For such $\lambda$ formula  \eqref{u_formula} is incorrect due to the small denominators. Therefore one have to study complex $\lambda$ and later pass to the limit $\operatorname{Im} \lambda \rightarrow 0$, see details, e.g., in Lecture 6 of \cite{M94} or in \cite{aizenman2015random}.

Our goal is to understand, under what conditions there exists an isolated positive eigenvalue $\lambda_0$ of the operator $H_\mu$, perturbed by the reproduction potential $(\Lambda - \mu)\delta_0(x)$. In such case \eqref{u_formula} is well defined and Lemma \ref{lemma_4} is probably new. Let us stress that $u(x)$ is the resolvent kernel of $H_\mu$ with some normalization. In fact $u(x) = \frac{R_\lambda(0,x)}{R_\lambda(0,0)}$ for $\lambda > 0$.
\end{Remark}

\begin{Remark}
The essential spectrum is the non-random support of the random spectral measure of $H_\mu(\omega)$. Under condition that r.v. $\mu(x)$ has absolutely continuous distribution, it follows from the general theory of one-dimensional random Schrödinger operator on $l_2(\mathbb Z)$ that the spectral measure is pure point and eigenfunctions are almost surely exponentially decreasing (exponential localization).

This result is very old, see details in \cite{M94,aizenman2015random,kunz1980spectre,carmona2012spectral}, and we will not discuss this topic. Our case is the analysis of the spectral bifurcation: existence and non-existence of the positive eigenvalue.
\end{Remark}

\section{Condition for almost certainly supercritical BRW behavior} \label{sec5}
Let us calculate the environments-independent interval in which the eigenvalue of the problem \eqref{eq3.1} lies. For this purpose, let us consider the most ``good'' and the most ``bad'' realizations of the environments. Namely, let us put $\mu = 0$ at all points, then let us put $\mu = c$ at all points.

\begin{Theorem}\label{lemma_6}
The value $P(\Lambda,\varkappa,F_{\mu})$ is equal to one if and only if 
the following condition is satisfied for the parameters of the BRW 
    \begin{equation}
    \label{lemma_6_cond}
        \Lambda \geqslant \sqrt{ (\varkappa + c)^2 - \varkappa ^2} - c.
    \end{equation}
If the condition \eqref{lemma_6_cond} is satisfied, then for any realization of the environments $\omega$, the eigenvalue of $\lambda(\omega)$ lies in the interval 
   $$\lambda(\omega) \in \bigl[\sqrt{(\Lambda + c)^2 + 
\varkappa^2} - (\varkappa + c); \sqrt{\Lambda^2 + \varkappa^2} - \varkappa 
\bigr].$$
\end{Theorem}
\textit{Proof}. Consider equation \eqref{eq3.2}:
\begin{equation}
\begin{aligned}
& \left( \varkappa \Delta + \Lambda - \lambda \right)u(0) = 0, \\
& u(0) =1.
\end{aligned}
\end{equation}
Given the notation $E = \varkappa + \lambda$, it can be rewritten as follows:
\begin{equation}
\label{brooklyn_baby}
\left(\varkappa \bar\Delta u \right) (0) + \Lambda - E = 0. \\
\end{equation}
In dimension $d = 1$, the expression \eqref{brooklyn_baby} takes a simpler form:
\[
\frac {\varkappa} {2} (u(1) + u(-1)) + \Lambda = E. \\
\]
Moreover, if $\mu$ is equal to some constant at all points, then $u(1) = u(-1)$ and the expression is further simplified:
\begin{equation}
\label{eigen_func_d_1}
\varkappa u(1) + \Lambda = E. \\
\end{equation}
Note that for arbitrary environment $\omega$, the solution $u(1)$ is bounded above and below by the solutions for the environments in which $\mu \equiv 0$ and $\mu \equiv c$. Let us find these estimates.

Let $\mu(x,\omega)$ be equal to some constant $c_1$ at all points. In this case, $u(1)$ is defined using equation \eqref{u_formula} as follows:
 \begin{equation}
 \label{u_sample_d1}
        u(1) = \sum\limits_{\gamma: 1 \rightarrow 0 } \prod\limits_{z \in \gamma} \left( \frac{\varkappa / 2}{\mu(z,\omega) + E} \right) = \sum\limits_{\gamma: 1 \rightarrow 0 } \prod\limits_{z \in \gamma} \left( \frac{\varkappa / 2}{E + c_1} \right) = \sum\limits_{\gamma: 1 \rightarrow 0 } \prod\limits_{z \in \gamma} \left( \frac{\varkappa / 2}{E_1} \right),
    \end{equation} 
where $E_1=E + c_1$.

Reasoning using the reflection principle as in the proof of Lemma~\ref{lemma_4} (see Subsection~\ref{sec9_2}) allows us to write out the series in the expression \eqref{u_sample_d1} exactly. First, let us compute $L(1,0,n)$, that is the number of paths that start at $1$, end at $0$, contain $n$ points, and do not intersect $0$. Note that $L(1,0,1) = 1$, and for the remaining odd $n$ according to reasoning \eqref{d_1_number_of_paths} in Subsection~\ref{sec9_2} the following is true:
$$L(1,0,n) = C_{n-1}^{\frac{n-1}{2}}-C_{n-1}^{\frac{n+1}{2}} =\frac{2}{n+1} C_{n-1}^{\frac{n-1}{2}}, \quad n = 3,5,\dots.$$
 
Thus, let us write out $u(1)$ from the expression \eqref{u_sample_d1}:
 \begin{multline}
 \label{eq6.0}
        u(1) = \sum\limits_{\gamma: 1 \rightarrow 0 } \prod\limits_{z \in \gamma} \left( \frac{\varkappa / 2}{ E_1} \right) = 
       \sum_{n = 1,3\dots}  L(1,0,n) \cdot \left(\frac{\varkappa/2}{E_1}\right)^n \\ 
        =\sum_{n = 1,3\dots} \frac{2}{n+1}  C_{n-1}^{\frac{n-1}{2}}  \cdot \left(\frac{\varkappa/2}{E_1}\right)^n 
        = 
        \sum_{m = 0,2\dots} \frac{2}{m+2}  C_{m}^{m/2}  \cdot \left(\frac{\varkappa/2}{E_1}\right)^{m+1} \\
            =\frac{\varkappa/2}{E_1} \sum_{k = 0}^\infty \left( \frac{C_{2k}^{k}}{1+k} \right)\cdot \left(\frac{\varkappa/2}{E_1}\right)^{2k}. 
    \end{multline}
For convenience, we denote $\frac{\varkappa/2}{E_1}$ by $a$. The coefficient $\frac{C_{2k}^{k}}{1+k}$ is the Catalan number, so the series \eqref{eq6.0} can be calculated exactly, see e.g. \cite{catalan}:
   \begin{multline}
    \label{eq6.0.1}
        u(1) = a \sum_{k = 0}^\infty \left( \frac{C_{2k}^{k}}{1+k} \right)\cdot a^{2k} = a \frac{1-\sqrt{1- 4a^2}}{2a^2} =\\ \frac{1}{2a} - \sqrt{\frac{1}{4a^2}-1} = \frac{E_1}{\varkappa} - \sqrt{ \left(\frac{E_1}{\varkappa}\right)^2-1}. 
    \end{multline}

Thus, the expression \eqref{eigen_func_d_1} takes the following form:
    \begin{equation*}
        E = \Lambda + {E_1} - \sqrt{ {E_1}^2-\varkappa^2}.
    \end{equation*}
We use the notation $E$, and the expression takes the form:
 \begin{equation*}
        E = \Lambda + (E + c_1) - \sqrt{(E + c_1)^2-\varkappa^2} .
    \end{equation*}
From here we can calculate that
$$ \Lambda + c_1 = \sqrt{ (E + c_1 ) ^2 - \varkappa^2}$$
or, finally,
\[
        \lambda = \sqrt{ (\Lambda + c_1) ^2 + \varkappa ^2} - (c_1 + \varkappa).
\]

Substituting $c_1 = 0$ and $c_1 = c$ completes the proof of the lemma. \qed

\section{Upper estimate for $P(\Lambda,\varkappa,F_{\mu})$.}  \label{sec6}

 In the previous section, we found that under the condition $ \Lambda \geqslant \sqrt{ (\varkappa + c)^2 - \varkappa ^2} - c$ the BRW is a.s. supercritical, i.e. $P(\Lambda,\varkappa,F_{\mu}) =1$. The goal of this and the next section is to give estimates for $P(\Lambda,\varkappa,F_{\mu})$ when this condition is violated. 
 
 To obtain an estimate from above, let us fix a non-random ``poor'' environment, and find out when it does not generate a positive eigenvalue. The poorer environments also do not generate an eigenvalue. If the probability of generating family of poor environments is $P_1$, then $P(\Lambda,\varkappa,F_{\mu})<1-P_1$. 
 In this paper, we consider the simplest case: a environment that takes some negative values at points neighboring zero.
\begin{Lemma} \label{lemma_9}
Consider a environment $\omega_1$  in which points neighboring from zero have killing intensities equal to $\mu_{1}$ and $\mu_{-1}$. A positive eigenvalue in this environments exists if and only if the following condition is met:
\begin{equation}
\label{eq9.1_main}
    \Lambda > \frac{\mu_1 + \mu_{-1} + 2\sigma \mu_1\mu_{-1}}{(1+\sigma \mu_1)(1+\sigma \mu_{-1})},
\end{equation}
where $\sigma = \frac{1}{\varkappa/2}$ for $z \in \mathbb{R}$.
\end{Lemma}
Let us give the general idea of the proof. The eigenvalue problem for the considered environment has the following form: 
\begin{equation}
\label{eq8.1_main}
\begin{aligned}
& \left( \varkappa \Delta +V(x, \omega) \right)u(x) = \lambda u(x), \\
& u(0) =1,
\end{aligned}
\end{equation}
where 
\begin{equation*}
V(x,\omega) = 
    \begin{cases}
         \Lambda,& \quad x = 0; \\
         -\mu_{1},& \quad x = 1; \\
         -\mu_{-1},&  \quad x = -1; \\
         0,& \quad |x| \geqslant 2.
    \end{cases}
\end{equation*}
In the appendix we show that $\psi(x)$ has the following form:
\begin{equation}
\label{eq8.4_main}
\psi(x) = 
\begin{cases}
    1,& \quad x =  0; \\
    C_1e^{-kx},& \quad x \geqslant 1; \\
    C_{-1}e^{kx},& \quad x \leqslant 1,
\end{cases}
\end{equation}
where $C_{\pm1}$ and $k$ are some positive constants. \\

Then we substitute \eqref{eq8.4_main} into \eqref{eq8.1_main} and derive the condition that is equivalent to the existence of positive eigenvalue $\lambda$ with corresponding eigenfunction $\psi_{\lambda}$. It turns out that this is the condition \eqref{eq8.1_main}, which completes the proof of the lemma. The proof is quite technical and is given in Subsection~\ref{sec9_3}. 

Now consider the general set of environments $\Omega_1 = \left\{\omega \in \Omega : \mu(1,\omega) = 
\mu_{1}, \mu(-1,\omega) = \mu_{-1} \right\}$. Note that the average 
number of particles in the non-random environment $\omega_1$ is a.s. greater than the average number of particles 
population in any environment from $\Omega_1$. Suppose that the condition \eqref{eq9.1_main} is satisfied for $\omega_1$. In such case, nothing can be said about the eigenvalues of the environments from $\Omega_1$. Suppose that the condition \eqref{eq9.1_main} is not satisfied for $\omega_1$. 
Then, according to the previous Lemma~\ref{lemma_9}, there is no positive eigenvalue for $\omega_1$ and hence there is no positive eigenvalue for all environments from $\Omega_1$. 

Let us denote the event ``condition \eqref{eq9.1_main} is met'' by $A$ and write the previous reasoning more formally:
\begin{multline}
\mathbb{P} \{\exists \ \lambda(\omega) >0 \} = \mathbb{P} \{\exists \ \lambda(\omega) >0 | A\}\mathbb{P}(A) + 
\mathbb{P} \{\exists \ \lambda(\omega) >0 |\bar A\}\mathbb{P}(\bar A)\\ = \mathbb{P} \{\exists \ \lambda(\omega) >0 | A\}\mathbb{P}(A) + 0 \cdot \mathbb{P}(\bar A) \leqslant \mathbb{P}(A).
\end{multline}

The event ``condition \eqref{eq9.1_main} is met'' for a random environment is written as follows:
$$\mathbb P (A) = \mathbb P\left\{\Lambda > \frac{\xi_1 + \xi_2 + 2\sigma \xi_1\xi_{2}}{(1+\sigma 
\xi_1)(1+\sigma \xi_{2})} \right\},$$ where $\xi_i$ are independent copies 
$\mu(x,\omega)$.
Thus, we obtain the following theorem.

\begin{Theorem} \label{Theorem_2}
The following upper bound estimate is true:
   $$P(\Lambda,\varkappa,F_{\mu}) \leqslant \mathbb{P} 
\left\{\Lambda > \frac{\xi_1 + \xi_2 + 2\sigma \xi_1\xi_{2}}{(1+\sigma 
\xi_1)(1+\sigma \xi_{2})} \right\},$$ where $\xi_i$ are independent copies 
$\mu(x,\omega)$.
\end{Theorem} 

\section{Lower estimate for $P(\Lambda,\varkappa,F_{\mu})$.}  \label{sec7}

The first way for obtaining a lower estimate for $P(\Lambda,\varkappa,F_{\mu})$ is to consider some convenient function $\psi(x)$ and examine the quadratic form $\left( H(\omega)\psi,\psi \right)$. If for some $a>0$ the quadratic form $\left( H(\omega)\psi,\psi \right)$ is positive with probability $p_a$, then the operator $H(\omega)$ has a positive eigenvalue with probability $p_a$ at least. We have taken the simple function $\psi(x) = 2^{-a|x|}, x\in \mathbb Z$ and this reasoning leads to the following theorem.

\begin{Theorem} \label{Theorem_3} The following estimate from below is true:
$$P(\Lambda,\varkappa,F_{\mu}) \geqslant \max\limits_{a\in (0;\infty)} \mathbb P \left( \omega: \Lambda > \varkappa\frac{(2^a-1)}{(2^a+1)}  + \sum\limits_{\substack{x = -\infty; \\ \ x\not= 0}}^{\infty} \frac {\mu(x,\omega)}{4^{a|x|}} \right).$$
In particular, for $a=1:$
$$P(\Lambda,\varkappa,F_{\mu}) \geqslant \mathbb P \left( \omega: \Lambda > \frac{\varkappa}{3}  + \sum\limits_{\substack{x = -\infty; \\ \ x\not= 0}}^{\infty} \frac {\mu(x,\omega)}{4^{|x|}} \right).$$
\end{Theorem}

The proof of the theorem requires direct investigation of the quadratic form $\left( H(\omega)\psi,\psi \right)$ for the function $\psi(x) = 2^{-a|x|}, x\in \mathbb Z$, which is a technical task, and so it is placed in Subsection~\ref{sec9_4}.

The second way to obtain an upper estimate of $P(\Lambda,\varkappa,F_{\mu})$ uses the idea of Lemma~\ref{lemma_9}. We consider a non-random killing environment of simplified form that can form <<islands>> around zero without killing. For this environment, we study the eigenvalue problem, and then generalize the conclusion to all environments that are ``better'' than the one under consideration.

First, let us denote $P(\mu(x,\omega)) = 0$ by $p$. Random variables $\mu(x,\omega)$ can form an ``island'' around zero with probability $p^{2l}$. Let us denote such case by $\Omega_l$:
$$\Omega_l = \left\{\omega \in \Omega : \mu(i,\omega) = 0, \forall i\in -l,\dots, l \right\}.$$
Let us use an idea from Lemma~\ref{lemma_9} and consider a non-random environment $\omega_l$ of the following form:
\begin{equation*}
    \mu(x,\omega_l) = \begin{cases}
0& \text{for } x \in -l,\dots, l; \\
c& \text{for } x \notin -l, \dots, l.
\end{cases}
\end{equation*}
The environment $\omega_l$ admits a direct calculation of the condition on the positivity of the eigenvalue of the corresponding operator, which is presented by the following lemma. The proof of the lemma is technical and is therefore placed in Subsection~\ref{sec9_5}.

\begin{Lemma} \label{add_1}
If a positive eigenvalue exists for all $\omega \in \Omega_l$, then it is bounded from below by a solution with respect to $\lambda$ of the following equation:
\begin{equation}
\label{eq10_1}
  \frac{2\alpha \varkappa }{1+\sqrt{1-4\alpha^2}} + \varkappa \alpha^{2l}  \cdot R(\alpha,\beta) + \Lambda -\varkappa - \lambda = 0, 
\end{equation}
where $\alpha = \frac{\varkappa/2}{\varkappa + \lambda}$, $\beta = \frac{\varkappa/2}{c+\varkappa + \lambda}$, and the expression $R$ is defined as follows:
\begin{equation}\label{eq_a1_1__}
 R(\alpha,\beta) = \sum\limits_{k=0}^\infty \left(\beta^{2k+1} - \alpha^{2k+1} \right) C_{k+l},
\end{equation}
where $ C_{n}$ denotes the $n$-th Catalan number.
If the series in equation \eqref{eq_a1_1__} does not converge  then there exists a $\omega \in \Omega_l$ for which there is no positive eigenvalue.
\end{Lemma}

Now, using Lemma~\ref{add_1} we find the smallest positive number $\hat l$ such that equation \eqref{eq10_1} admits a positive solution. By the lemma, all environments of $\Omega_{\hat l}$ will have positive eigenvalues. Therefore the probability $P(\Lambda,\varkappa,F_{\mu})$ is at least equal to the probability of generating an environment from $\Omega_{\hat l}$ or, equivalently, the probability of generating a $\hat l$-island. Finally, the probability of generating a $\hat l$-island is $p^{2 \hat l}=\left(\mathbb P\{\mu(x,\omega) = 0\} \right)^{2 \hat l}$, which leads to the following theorem.

\begin{Theorem} \label{Theorem_4}
There is the following estimation from below:
\begin{equation*}
P(\Lambda,\varkappa,F_{\mu}) \geqslant \left(\mathbb P\{\mu(x,\omega) = 0\} \right)^{2\hat l},
\end{equation*}
where $\hat l \in \mathbb{N}$ is the smallest number for which the expression described below admits a positive solution. If there is no such $\hat l$, then $P(\Lambda,\varkappa,F_{\mu}) = 0.$
\begin{equation}
  \frac{2\alpha \varkappa }{1+\sqrt{1-4\alpha^2}} + \varkappa \alpha^{2l}  \cdot R(\alpha,\beta) + \Lambda -\varkappa - \lambda = 0, 
\end{equation}
where $\alpha = \frac{\varkappa/2}{\varkappa + \lambda}$, $\beta = \frac{\varkappa/2}{c+\varkappa + \lambda}$, and the expression $R$ is defined as follows:
\begin{equation}\label{eq_a1_1}
 R(\alpha,\beta) = \sum\limits_{k=0}^\infty \left(\beta^{2k+1} - \alpha^{2k+1} \right) C_{k+l},
\end{equation}
where $ C_{n}$ denotes the Catalan number.
\end{Theorem}

At first sight, Theorem~\ref{Theorem_4} is useless due to its excessive complexity. However, unlike Theorem~\ref{Theorem_3}, it offers a concrete numerical algorithm for estimating $P(\Lambda,\varkappa,F_{\mu})$ based on non-Monte Carlo method. Moreover, these algorithm will run fast because of the exponentially fast convergence of the series used in the theorem.

\section{Conclusion}

In this paper we propose approaches that can be used for analysis of BRW in a random environment. We derived the condition of almost sure supercriticality of the BRW. In case this condition is violated, we obtain lower and upper estimates for the probability of supercriticality. The estimates can be improved: in Theorem~\ref{Theorem_2} one can consider not two but several points in the neighborhood of zero, in Theorem~\ref{Theorem_3} one can consider a function of more general form, and in Theorem~\ref{Theorem_4} one can consider an island consisting of points with small but non-zero positive killing intensity.

\section{Proofs}\label{sec9}

\subsection{Continuation of the proof of Lemma~\ref{lemma3}} {\label{sec9_1}}
Let us recall, we need to prove the inclusion $\sigma(H_\mu) \supseteq [-2\varkappa-c;0]. $ Let us use Lemma~\ref{lemma1}, and for each $\lambda \in [-2\varkappa-c;0]$ construct a sequence of ``almost eigenfunctions'' $\{ f_n(x) \},$  $f_i(x) \in L_2 (\mathbb Z ^d)$. Note that we construct a sequence function for each fixed $\omega$, i.e., $\{ f_n(x) \}$ = $\{ f_n(x,\omega) \}$. 

Let us represent $\lambda$ as $\lambda = a + b$, $a \in [-2\varkappa;0]$ $b \in [-c;0]$. Let us construct $f_n$ such that in some sense $\varkappa \Delta f_n \approx a f_n$ and simultaneously $-\mu(x,\omega) f_n \approx b f_n$. Then 
$$\left(\varkappa \Delta -\mu(x,\omega) \right) f_n \approx (a + b) f_n = \lambda f_n.$$

The condition $\varkappa \Delta f_n \approx a f_n$ requires that the function be ``almost everywhere'' similar to $\exp \{ i(\overrightarrow{\phi}, x)\}$ of Lemma~\ref{lemma2} with a suitable $\overrightarrow{\phi}$. The condition $-\mu(x,\omega) f_n \approx b f_n$ requires that the function be nonzero only on the region where $-\mu(x,\omega) \approx b$.

A candidate function satisfying both conditions looks like this:
\[
 f_n(x) = f_n(x,\omega) = \frac{1}{\sqrt{|B_n(\omega)|}} \exp \left\{ i(\overrightarrow{\phi}, x) \right\} I\{B_n(\omega)\},    
\]
where the random set $B_n(\omega) = B_n$ contains the points $x \in \mathbb{Z}^d$, such that $-\mu(x,\omega) \in \left[b-\frac 1n, b + \frac 1n \right],$ and the multiplier ${1}/{\sqrt{|B_n|}}$ is needed to normalize the function.

Lemma~\ref{lemma1}  additionally requires orthogonality of almost eigenfunctions. Hence it should be required that $B_m \cap B_n = 0$ for $n \not= m$. Furthermore, it will turn out in the proof that it should be required in advance that $|B_n| \rightarrow \infty$. Let us show the existence of the required sets $\{ B_n \}$. 

Let us fix an arbitrary number $n$. Recall that the density of $-\mu(x.\omega)$ is positive on the interval $[-c;0]$. According to the Borel-Cantelli lemma, for an arbitrary realization of $\omega$ there exists a system of non-intersecting balls $\{C_i(n)\}_{i=1}^\infty$ consisting of lattice points $x \in \mathbb{Z}^d$ such that:
\begin{equation}
\begin{aligned}
  x \in C_i(n)&\Rightarrow -\mu(x,\omega) \in \left[b-\frac 1n, b + \frac 1n \right], \\
  |C_i(n)| &\rightarrow \infty\quad\text{for } i \rightarrow \infty.
\end{aligned}
\end{equation}

Now the system of sets $\{ B_n \}$ can be constructed by induction. Let the system $\{ B_n \}$ be constructed up to the number $k$. Let us construct the system $\{C_i(k+1)\}$ described above. As the set $B_{k+1}$, we take any set of $\subset \{C_i(k+1)\}$ that is farther from zero than all points from $B_1, \dots B_n$. Thus, the induction is complete and the system $\{ B_n \}$ is constructed.

Let us verify that the functions $\{f_n\}$ are almost eigenfunctions. First, consider the action of the operator $\varkappa\Delta$ on the function $f_n$:
$$
    \varkappa\Delta f_n = \frac{1}{\sqrt{|B_n|}}  \varkappa\Delta \exp \{ i(\overrightarrow{\phi}, x)\} I\{B_n\}.
$$
If the points $x-1, x, x, x+1$ lie inside $B_n$, then the mixing operator acts on its eigenfunction:
$$
\varkappa\Delta f_n (x) = a \frac{1}{\sqrt{|B_n|}} \exp \{ i(\overrightarrow{\phi}, x) \} = a f_n (x).
$$
If all points $x-1, x, x+1$ lie outside $B_n$, then the mixing operator acts on the null function and also $\varkappa\Delta f_n (x) = 0 = af_n (x).$

Let at least one of their points $x-1, x, x+1$ lie on the boundary of $B_n$. In this case, there remains a non-zero function $f_n^{res}$ after applying the operator:
\begin{equation} \label{first_approx}
    \varkappa\Delta f_n (x) = a f_n (x) + f_n^{res}(x).
\end{equation}
The function $f_n^{res}$ reflects the ``error'' of the operator on the boundary of $B_n$ with respect to the  operator multiplying by $a$. This function is non-zero only at a finite number of points $C_d$, depending on the dimensionality but not on n. The norm $f_n^{res}$ is bounded from above by $ C_d/\sqrt{|B_n|}$, which tends to zero when $n \rightarrow \infty$.

Now consider the action of the operator $-\mu(x, \omega)$ on the function $f_n$. On the region $\{ B_n \}$, the function $ \mu(x, \omega)$ takes values in the interval $\left[b-\frac 1n, b + \frac 1n \right]$, so:
\begin{multline} \label{second_approx}
        -\mu(x, \omega) f_n (x) = - \mu(x, \omega) \frac{1}{\sqrt{|B_n|}} \exp \{ i(\overrightarrow{\phi}, x)\} I\{B_n\} \\
    =b \frac{1}{\sqrt{|B_n|}}  \exp \{ i(\overrightarrow{\phi}, x)\} I\{B_n\} + g_n^{res}(x) = b f_n (x)+ g_n^{res}(x).
\end{multline}
The function $g_n^{res}$ reflects the <<error>> of the operator on the area $B_n$ with respect to the multiplication operator on $b$. The norm $g_n^{res}$ is bounded from above by the value $ 1/n$, which tends to zero when $n \rightarrow \infty$.

Putting together the expressions \eqref{first_approx} and \eqref{second_approx}, we get:
\begin{equation*}
 \left\|H_\mu f_n- \lambda f_n\right\| =   \left\| (a+b) f_n- (a+b) f_n + f_n^{res} + g_n^{res}\right\| \rightarrow 0, \qquad n \rightarrow \infty.
 \end{equation*}
Thus, $\{f_n\}$ is the desired sequence of almost eigenfunctions, and $\lambda \subset \sigma(H).$ The number $\lambda$ was taken arbitrarily from the interval $[-2\varkappa-c,0]$, hence:
\begin{equation}
\sigma(H_\mu) \supseteq [-2\varkappa-c;0]. \qed
\end{equation}

\subsection{Proof of Lemma~\ref{lemma_4}} {\label{sec9_2}}
Let us recall the formulation of Lemma~\ref{lemma_4}: \\
        The solution to equation \eqref{eq3.4} when $x \not = 0$ is given by the following formula:
    \begin{equation}
    \label{u_formula_add}
        u(x) = \sum\limits_{\gamma: x \leadsto 0 } \prod\limits_{z \in \gamma} \left( \frac{\varkappa / 2d}{\mu(z,\omega) + E} \right),
    \end{equation}
    where by $\gamma: a \leadsto b = \{ a = x_1 ,\dots x_{n} \not= b\}$ denotes the path from point $a$ to point $b$ through the neighboring points of the lattice, and: a) the path does not intersect $0$; b) point $b$ is not considered to be included in the path $\gamma$.    
    The solution given by formula \eqref{u_formula_add} makes sense for any $\lambda >0$ in any dimension $d \in \mathbb{N}$.

 Note that, for a path $\gamma: x\leadsto 0$, the symbol $|\gamma|$ denotes the length of the path in the sense of ``number of points in $\gamma$ excluding zero'' or ``number of steps from $x$ to $0$'' which is the same. For simplicity of notation, let us prove the lemma for the case of a one-dimensional lattice $d = 1$. Let us study the action of the operator $\varkappa \bar \Delta$ on the function u(x) when $x \not= 0$:
\begin{equation}
\label{barkappa_action_add}
    \varkappa \bar \Delta u(x) = \frac{\varkappa}{2} \left( u(x+1) + u(x-1) \right),
\end{equation}
Note that the set of paths $\gamma: x \leadsto 0$ included in $u(x)$ decomposes into two subsets: paths $\gamma_+: x \leadsto x+1 \leadsto 0$ and paths $\gamma_-: x \leadsto x-1 \leadsto 0$. Thus,
\begin{multline}
\label{eq_u_long_add}
    u(x) = \sum\limits_{\gamma} (\cdot)= \frac{\varkappa / 2}{\mu(x,\omega) + E}  \sum\limits_{\gamma_+} (\cdot) + \frac{\varkappa / 2}{\mu(x,\omega) + E}  \sum\limits_{\gamma_-} (\cdot) \\ =\frac{\varkappa / 2}{\mu(x,\omega) + E} \left( u(x+1) + u(x-1) \right).
\end{multline}
Or, which is the same thing:
\begin{equation}
\label{u_representation_add}
    u(x+1) + u(x-1) = \frac{\mu(x,\omega) + E}{\varkappa / 2} u(x).
\end{equation} 
Combining \eqref{barkappa_action_add} and \eqref{u_representation_add} we get:
\begin{equation}\label{eq4_1_add}
    \varkappa \bar \Delta u(x) = \frac{\varkappa}{2} \left( u(x+1) + u(x-1) \right) = (\mu(x,\omega) + E) u(x).
\end{equation}
By virtue of \eqref{eq4_1_add}, the proof of the lemma in the one-dimensional case is complete:
\begin{equation*}
   \left( \varkappa \bar \Delta - ( \mu(x, \omega) + E) \right)u(x) = (\mu(x,\omega) + E) u(x) - (\mu(x, \omega) + E)u(x) =0.
\end{equation*}
In the multidimensional case, the reasoning remains exactly the same, except that the expression \eqref{eq_u_long_add} will contain paths on all lattice points neighboring $x$.

Now let us show the correctness of \eqref{u_formula_add} for any $\lambda>0$. For simplicity, we first consider the one-dimensional case $d = 1$.  We investigate the convergence of the series \eqref{u_formula_add}. Note that the following upper bound estimate is true and achievable: 
 \begin{equation}
 \label{upper_est}
        u(x) = \sum\limits_{\gamma: x \leadsto 0 } \prod\limits_{z \in \gamma} \left( \frac{\varkappa / 2}{\mu(z,\omega) + E} \right) \leqslant \sum\limits_{\gamma: x \leadsto 0 } \prod\limits_{z \in \gamma} \left( \frac{\varkappa }{2E} \right) = \sum\limits_{{\substack{\gamma: x \leadsto 0, \\ |\gamma| = n}}} \left( \frac{\varkappa}{2E} \right)^{n}L(x,0,n),
    \end{equation}
    where $L(x,0,n)$ is the number of paths of the form $x \leadsto 0$ that contain $n$ points. Note that if the parity of $x$ and $n$ does not coincide, then $L(x,0,n)$ converges to zero.

    Without restricting generality, let us assume $x>0$. Finding $L(a,b,k)$ is a standard problem for applying the reflection principle to discrete random walks. The answer is as follows:
    \[
     L(a,b,k) = C_k^{\frac{k+b-a}{2}}-C_k^{\frac{k+b+a}{2}}, \qquad a,b,n>0.
     \]
    Therefore:   
      \begin{equation} \label{d_1_number_of_paths}
        L(x,0,n) = L(x,1,n-1) = C_{n-1}^{\frac{n-x}{2}} -  C_{n-1}^{\frac{n+x}{2}}  = c_1 x^{c_2}2^{n}, \qquad n \rightarrow \infty,
    \end{equation} 
    where $c_1$ and $c_2$ are positive constants.
    
Thus, the series in the \eqref{upper_est} inequalities are geometric series:
  \begin{equation} \label{eq5_1}
          \sum\limits_{{\substack{\gamma: x \leadsto 0, \\ |\gamma| = n}}} \left( \frac{\varkappa}{2E} \right)^{n}L(x,0,n) < c_3 + c_4\sum\limits_{i = 1}^n  \left( \frac{\varkappa}{2E} \right)^{n} \cdot 2^{n}, 
   \end{equation}
    where $c_3$ and $c_4$ are positive constants.
    
    The series \eqref{eq5_1} converges when $ {\varkappa}/{E} < 1.$ Which, given the notation $E = \varkappa + \lambda$, can be rewritten as follows:
    \[
         \lambda > 0.    
    \]

In the case $d > 1$, the estimaion of \eqref{upper_est} takes the following form:
 \begin{equation}
 \label{upper_est_d2}
        u(x) = \sum\limits_{\gamma: x \leadsto 0 } \prod\limits_{z \in \gamma} \left( \frac{\varkappa / 2d}{\mu(z,\omega) + E} \right) \leqslant \sum\limits_{\gamma: x \leadsto 0 } \prod\limits_{z \in \gamma} \left( \frac{\varkappa }{2dE} \right) = \sum\limits_{{\substack{\gamma: x \leadsto 0, \\ |\gamma| = n}}} \left( \frac{\varkappa}{2dE} \right)^{n}L(x,0,n),
    \end{equation}
    where again $C_n$ is the number of $\gamma$ paths of length $n$, where the length is counted with respect to zero and the conditions from Lemma~\ref{lemma_4} are imposed on the path.

   Let us consider $n\gg x$, since the convergence of the series \eqref{upper_est_d2} depends on them alone. Note that when $n \gg x$, the number of trajectories $L(x,0,n) \sim L(0,0,n)$, $n \rightarrow \infty$. We denote by $L_0(0,0,n)$ the number of trajectories starting and ending at zero without the condition of non-intersection of zero. The event of a trajectory crossing the zero point in dimension $d>1$ is rare, so $L(0,0,n) \sim L_0(0,0,n)$ , $n \rightarrow \infty$.

   Let us fix $d$ movements ``up'' along each of the coordinates. Each path in $L_0(0,0,n)$ is defined by only $n/2$ steps, each of which can have one of the coordinate movements, i.e.:
   \begin{equation*}
       L_0(0,0,n) = d^{n/2}C_n^{n/2} \sim \left(2\sqrt d \right)^n, \quad n \rightarrow \infty. 
   \end{equation*}

   Proceeding as in the one-dimensional case, we obtain that the series in the estimaion of \eqref{upper_est_d2} is a geometric series:
     \begin{equation}\label{eq5_2}
      \sum\limits_{{\substack{\gamma: x \leadsto 0, \\ |\gamma| = n}}} \left( \frac{\varkappa}{2dE} \right)^{n}L(x,0,n) < c_5 + \sum\limits_{i = 1}^n  \left( \frac{\varkappa}{2dE} \right)^{n} \cdot \left(2\sqrt d \right)^n,
        \end{equation}
        where $c_5$ is a positive constant.
    The series \eqref{eq5_2} converges when $ {\varkappa}/{\sqrt d E} < 1.$ Which given the notation $E = \varkappa + \lambda$ can be rewritten as follows:
     \[
         \lambda > \varkappa \left(\frac 1 {\sqrt{d}} - 1 \right).    
    \]
    Therefore for $\lambda >0$ the series converges which completes the proof of the lemma. \qed

\subsection{Proof of Lemma~\ref{lemma_9}} \label{sec9_3}
Let us recall the formulation of Lemma~\ref{lemma_9}: \\
Consider a  environment $\omega_1$ in which points neighboring from zero have killing intensities equal to $\mu_{1}$ and $\mu_{-1}$. A positive eigenvalue in this environment exists if and only if:
\[
    \Lambda > \frac{\mu_1 + \mu_{-1} + 2\sigma \mu_1\mu_{-1}}{(1+\sigma \mu_1)(1+\sigma \mu_{-1})},
\]
where $\sigma = \frac{1}{\varkappa/2}$ for $z \in \mathbb{R}$.

 The eigenvalue problem for the considered environment is as follows: 
\begin{equation}
\label{eq8.1}
\begin{aligned}
& \left( \varkappa \Delta +V(x, \omega) \right)u(x) = \lambda u(x), \\
& u(0) =1,
\end{aligned}
\end{equation}
where 
\begin{equation*}
V(x,\omega) = 
    \begin{cases}
         \Lambda,& \quad x = 0; \\
         -\mu_{1},& \quad x = 1; \\
         -\mu_{-1},& \quad x = -1; \\
         0,& \quad |x| \geqslant 2.
    \end{cases}
\end{equation*}

First, let us show how the eigenfunction for this problem looks like in general form. We will use the forward and inverse discrete Fourier transforms, see, e.g., \cite{yarovaya2007}. The Fourier transform of the function $f$ is defined as follows:
$$\tilde f (\theta) = \sum\limits_{x \in \mathbb{Z}} e^{i\theta x} f(x).$$
The inverse Fourier transform is defined as follows:
$$ f (x) = \frac{1}{2\pi} \int\limits_{-\pi}^{\pi}\tilde f(\theta)e^{-i\theta x} d\theta.$$

Let us write for the operator $H$ the eigenvalue problem $\lambda$ with the corresponding eigenfunction $u$:

\begin{equation} \label{prev_Fourier}
    \varkappa\Delta u(x) + \Lambda \delta_0(x) u(x) - \mu_{-1}  \delta_{-1}(x) u(x) - \mu_{1}  \delta_{1}(x) u(x) = \lambda u(x).
\end{equation}
After applying the Fourier transform, the expression \eqref{prev_Fourier} takes the form
\[
    \varkappa (\cos(\theta)-1) \tilde u(x) + \Lambda u(0) - \mu_{-1} u(-1) e^{-i \theta}- \mu_{1} u(1) e^{i \theta} - \lambda \tilde u(x).
\]
The Fourier transform of the eigenfunction $\tilde u(x)$ is as follows:
\[
    \tilde u(x) = \frac{\Lambda u(0) - \mu_{-1} u(-1) e^{-i \theta}- \mu_{1} u(1) e^{i \theta}}{\lambda + \varkappa - \varkappa \cos \theta},
\]
and, finally, the solution $u(x)$ can be represented as
\[
     u(x) = \frac{1}{2 \pi} \int\limits_{-\pi}^{\pi} \frac{\Lambda u(0) - \mu_{-1} u(-1) e^{-i \theta}- \mu_{1} u(1) e^{i \theta}}{\lambda + \varkappa - \varkappa \cos \theta} e^{-i\theta x} d\theta.
\]
Calculating here the integral for $x\geqslant 1$ we obtain:
\begin{equation} \label{u_postpost}
     u(x) = -\mu_{-1}u(-1) \frac{w^{x-1}}{r} +\Lambda u(0) \frac{w^{x}}{r} -\mu_{1}u(1) \frac{w^{x+1}}{r},
\end{equation}
where $r = \sqrt{\lambda(\lambda + 2\varkappa)}$ and $w = \frac{\lambda + \varkappa - r}{\varkappa}$.

The expression \eqref{u_postpost} can be rewritten in a more convenient form:
\begin{equation}
     u(x) = w^{x}\left(-\mu_{-1} u(-1)\frac{1}{wr} + \Lambda u(0)\frac{1}{r} -\mu_1 u(1)\frac{w}{r}\right) = B_1 \cdot e^{-kx},
\end{equation}
where $B_1 = \left(-\mu_{-1} u(-1) \frac{1}{wr} + \Lambda u(0) \frac{1}{r} -\mu_1 u(1) \frac{w}{r}\right)$ and $ e^{-k} = w$. Function $w^x$ decreases as $x$ tends to infinity, since $u \in L_2(\mathbb{Z})$, therefore $k>0$.

Let us do exactly the same for $x \leqslant 1$ and put $f(0) = 1$ by normalization. We obtain that the eigenfunction must have the following form:
\[
\psi(x) = 
\begin{cases}
    1,& \quad x = 0; \\
    C_1e^{-kx},& \quad x \geqslant 1; \\
    C_{-1}e^{kx},& \quad x \leqslant -1,
\end{cases}
\]
where $C_{\pm1}$ and $k$ are positive constants. Let us show that there is a positive eigenvalue for this function if and only if the lemma condition is satisfied.

Let us write the problem \eqref{eq8.1} for the point $x \in [2;\infty)$:
\begin{equation*}
    \frac{\varkappa}{2}\psi(x+1) + \frac{\varkappa}{2}\psi(x-1) - \varkappa \psi(x) = \lambda \psi(x);
\end{equation*}
from here we can calculate that
\begin{equation*}
   \lambda = \frac\varkappa 2 \left( e^{-k} + e^k - 2 \right) = \varkappa \left( \cosh k - 1 \right) = 2 \varkappa \sinh^2(k/2) = k^2 + O(k^4), \quad k \rightarrow 0.
\end{equation*}
In particular, when $\lambda \rightarrow 0+$, it follows from the condition $k > 0$ that $k \rightarrow 0+$, that is
\begin{equation}
\label{z_limiting}
     \lambda \rightarrow 0+ \quad \rightarrow \quad e^k \rightarrow 1+.
\end{equation}

Now let us write the problem \eqref{eq8.1} for the point $x = 1$:
\begin{equation*}
    \frac{\varkappa}{2}\psi(2) + \frac{\varkappa}{2}\psi(0) - \varkappa \psi(1) - \mu \psi(1) = \lambda \psi(1).
\end{equation*}
From here, we can calculate that
\begin{equation*}
   C_1 = \frac{1}{1+e^{-k}\frac{\mu_1}{\varkappa/2}}.
\end{equation*}
For simplicity we denote $\frac{1}{\varkappa/2}$ by $\sigma$. let us perform similar reasoning for $x = -1$, obtaining:
\begin{equation*}
   C_{\pm1} = \frac{1}{1+\sigma \mu_{\pm 1} e^{-|k|}}.
\end{equation*}
Finally, let us write the problem \eqref{eq8.1} for $x = 0$: 
\begin{equation*}
    \frac{\varkappa}{2}\psi(1) + \frac{\varkappa}{2}\psi(-1) - \varkappa \psi(0) + \lambda \psi(0) = \lambda \psi(0).
\end{equation*}
From here, we can calculate that
\begin{equation*}
   \left(C_{1} + C_{-1} - 1\right)e^{-k} + \sigma \Lambda - e^k = 0
\end{equation*}
or, finally,
\begin{equation*}
    e^{2k} - \sigma \Lambda e^k - \left(\frac{1}{1+\sigma \mu_{1} e^{-k}} +\frac{1}{1+\sigma \mu_{- 1} e^{-k}} -1 \right) = 0.
\end{equation*}

First, for simplicity, let us put $\mu = \mu_{-1} = \mu_1$, also denote $e^k$ by $z$ and get the following expression:
\[
    z^2 - \sigma \Lambda z - \frac{2}{1+\sigma \mu \frac{1}{z}} + 1 = 0
\]
or
\begin{equation}
\label{z_polynom}
    z^3 - z^2 \sigma (\Lambda - \mu) - z(\sigma^2 \Lambda \mu + 1) + \sigma \mu = 0.
\end{equation}

The current problem is to write out conditions on $\Lambda, \mu$, and $\sigma$ that guarantee the positivity of $\lambda$. Let us use the expression \eqref{z_limiting} and note that \eqref{z_polynom} is a smooth function with respect to $z$, so we can put $z=1$ to find the limit solution at $z \rightarrow 1+$:
\begin{equation*}
    1 - \sigma (\Lambda - \mu) - (\sigma^2 \Lambda \mu + 1) + \sigma \mu = 0 \Leftrightarrow \Lambda = \frac{2\mu}{1+\sigma\mu}.
\end{equation*}
The eigenvalue $\lambda > 0$ exists when this condition is violated to the ``supercritical side'', i.e.
\begin{equation*}
    \Lambda > \frac{2\mu}{1+\sigma\mu}.
\end{equation*}
We obtain the conditions of the lemma under the assumption of $\mu_1 = \mu_{-1}$.

In the case of unequal $\mu_1$ and $\mu_{-1}$, the same solution method gives the condition of the lemma:
\[
    \Lambda > \frac{\mu_1 + \mu_{-1} + 2\sigma \mu_1\mu_{-1}}{(1+\sigma \mu_1)(1+\sigma \mu_{-1})}. \quad \qed
\]

\subsection{Proof of Theorem~\ref{Theorem_3}}\label{sec9_4}
Consider the function $\psi(x) = 2^{-a|x|}$. Let us denote $\varphi(x) = \left( H(\omega)\psi \right)(x)$ and directly calculate the quadratic form $(\varphi,\psi)= \left(\varphi(x,\omega),\psi(x)\right)$. First, let us calculate the function $\varphi(x):$
\begin{multline} \label{eq9_1}
    \varphi(x,\omega) = \varkappa\Delta\psi(x) + \Lambda\delta_0(x)\psi(x) - (1 - \delta_0(x)) \mu(x,\omega)\psi(x) \\
   =\frac{  \varkappa}{2} (\psi(x+1)+\psi(x-1)+2\psi(x)) + \Lambda\delta_0(x)\psi(0) - (1 - \delta_0(x)) \mu(x,\omega)\psi(x).
\end{multline}
Let us substitute into \eqref{eq9_1} the expression for $\psi(x)$ and consider separately the points $0$ and $x>0$:
\begin{align}\label{eq9_2}
    \varphi(0,\omega)& =\frac{\varkappa}{2} (2^{-a} + 2^{-a} - 2) + \Lambda = \varkappa (2^{-a}-1)  + \Lambda;\\
\nonumber
    \varphi(x,\omega) &=\frac{\varkappa}{2} ( 2^{-a}\cdot2^{-ax}+ 2^{a} \cdot 2^{-ax} -2\cdot2^{-ax}) - \mu(x,\omega)2^{-ax}=\\
\label{eq9_3}     &\qquad\qquad\qquad\qquad\qquad\frac{\varkappa}{2} 2^{-ax}( 2^{-a}+2^{a}-2) - \mu(x,\omega)2^{-ax}.
\end{align}

Using \eqref{eq9_2} and \eqref{eq9_3} we calculate  the required quadratic form:
\begin{multline} \label{eq9_4}
   \left(\varphi(x,\omega),\psi(x)\right) = \sum\limits_{\substack{x = -\infty; \\ \ x\not= 0}}^{\infty} \varphi(x)\psi(x) + \varphi(0)\psi(0) \\ 
    =\sum\limits_{\substack{x = -\infty; \\ x\not= 0}}^{\infty} \left(\frac{\varkappa}{2} 2^{-a|x|}( 2^{-a}+2^{a}-2) - \mu(x,\omega)2^{-a|x|}\right)\cdot 2^{-a|x|} + \varkappa (2^{-a}-1)  + \Lambda \\ 
    =\frac{\varkappa}{2} ( 2^{-a}+2^{a}-2) \sum\limits_{\substack{x = -\infty; \\ \ x\not= 0}}^{\infty} 2^{-2a|x|} - \sum\limits_{\substack{x = -\infty; \\ x\not= 0}}^{\infty} \frac {\mu(x,\omega)}{2^{2a|x|}} + \varkappa (2^{-a}-1)  + \Lambda \\
    =-\varkappa \frac{2^{-a}-1}{(2^{-a}+1)2^a}+ \varkappa (2^{-a}-1)  + \Lambda  - \sum\limits_{\substack{x = -\infty; \\ x\not= 0}}^{\infty} \frac {\mu(x,\omega)}{2^{2a|x|}} \\ = -\varkappa\frac{(2^a-1)}{(2^a+1)}  + \Lambda  - \sum\limits_{\substack{x = -\infty; \\ x\not= 0}}^{\infty} \frac {\mu(x,\omega)}{2^{2a|x|}}.
\end{multline}

If $\left(\varphi(x,\omega),\psi(x)\right) > 0 $, then by virtue of Section~\ref{section_4}, the operator $H(\omega)$ has a positive eigenvalue. Given the expression \eqref{eq9_4}, the condition for the positivity of the quadratic form can be rewritten in the following form:
\begin{equation*}
   \Lambda > \varkappa\frac{(2^a-1)}{(2^a+1)} + \sum\limits_{\substack{x = -\infty; \\ \ x\not= 0}}^{\infty} \frac {\mu(x,\omega)}{2^{2a|x|}}.
\end{equation*}
By substituting $a = 1$ we get: 
\begin{equation*}
   \Lambda > \frac{\varkappa}{3} + \sum\limits_{\substack{x = -\infty; \\ \ \ x\not= 0}}^{\infty} \frac {\mu(x,\omega)}{4^{|x|}}. \qed
\end{equation*}

\subsection{Proof of Lemma~\ref{add_1}} \label{sec9_5}
Recall the formulation of the lemma:\\
Consider a set of $\Omega_l$ including environments that have $l$-islands around zero: 
$$\Omega_l = \left\{\omega \in \Omega : \mu(i,\omega) = 0,~ \forall i\in -l,\ldots,l \right\}.$$
If a positive eigenvalue exists for all $\omega \in \Omega_l$, then it is bounded from below by a solution with respect to $\lambda$ of the following equation:
\begin{equation*}
  \frac{2\alpha \varkappa }{1+\sqrt{1-4\alpha^2}} + \varkappa \alpha^{2l}  \cdot R(\alpha,\beta) + \Lambda -\varkappa - \lambda = 0, 
\end{equation*}
where $\alpha = \frac{\varkappa/2}{\varkappa + \lambda}$, $\beta = \frac{\varkappa/2}{c+\varkappa + \lambda}$, and the expression $R$ is defined as follows:
\begin{equation}\label{eq_a1_1_add}
 R(\alpha,\beta) = \sum\limits_{k=0}^\infty \left(\beta^{2k+1} - \alpha^{2k+1} \right) C_{k+l},
\end{equation}
where $ C_{n}$ denotes the Catalan number.
If the series in equation \eqref{eq_a1_1_add} does not converge then there exists a $\omega \in \Omega_l$ for which there is no positive eigenvalue.

\textit{Proof}. 
For convenience in the proof, let us rewrite formula \eqref{lemma_4} for $d=1$:
    \begin{equation}
    \label{u_formula_again}
        u(x) = \sum\limits_{\gamma: x \leadsto 0 } \prod\limits_{z \in \gamma} \left( \frac{\varkappa / 2}{\mu(z,\omega) + E} \right).
    \end{equation}

By virtue of equation \eqref{eigen_func_d_1}, we are interested in the quantity $u(1)$, for which the paths from $1$ to $0$ are important.
By the lemma condition, for any $\omega \in \Omega_l$ there exists an $l$-island around zero. Therefore, every trajectory from $1$ to $0$ of length less than $2l$ will not leave this island, and every trajectory of length greater than $2l$ must spend at least $2l$ steps in this island. We divide all trajectories into two families: trajectories of length less than or equal to $2l$ trajectories of length greater than $2l$. The contribution of each of the smaller trajectories to the sum \eqref{u_formula_again} is exactly $\left(\frac{\varkappa/2}{0+E}\right)^{|\gamma|}$. The contribution of each of the large trajectories to the sum \eqref{u_formula_again} is at least $\left(\frac{\varkappa/2}{0+E}\right)^{2l} \cdot \left(\frac{\varkappa/2}{c+E}\right)^{|\gamma|-2l}$. Thus,
\begin{multline}
\label{add_1_1}
    u(1) = \sum\limits_{\gamma: 1 \rightarrow 0 } \prod\limits_{z \in \gamma} \left( \frac{\varkappa / 2}{\mu(z,\omega_l) + E} \right)\\     \geqslant\sum\limits_{{\substack{\gamma: 1 \rightarrow 0, \\ |\gamma| < 2l}}}   \left( \frac{\varkappa / 2}{0 + E} \right)^{|\gamma|} L(0,1,|\gamma|) + \sum\limits_{{\substack{\gamma: 1 \rightarrow 0, \\ |\gamma| \geqslant 2l}}}  \left( \frac{\varkappa / 2}{0 + E} \right)^{2l} \left( \frac{\varkappa / 2}{c + E} \right)^{ 
    |\gamma| - 2l}L(0,1,|\gamma|).
\end{multline}
Let us rewrite equation \eqref{add_1_1} by introducing the notations $\alpha = \frac{\varkappa/2}{E}$, $\beta = \frac{\varkappa/2}{c+E}$:
\begin{equation}
\label{add_1_2}
    u(1) \geqslant
    \sum\limits_{{\substack{\gamma: 1 \rightarrow 0, \\ |\gamma| < 2l}}}    \alpha ^{|\gamma|} L(0,1,|\gamma|) + \sum\limits_{{\substack{\gamma: 1 \rightarrow 0, \\ |\gamma| \geqslant 2l}}}   \alpha^{2l} \beta^{ 
    |\gamma| - 2l}L(0,1,|\gamma|).
\end{equation}
Let us simplify equation \eqref{add_1_2} in the way \eqref{eq6.0} does, yielding the following:
\begin{multline}
\label{add_1_3}
    u(1) \geqslant \sum\limits_{k=0}^{l-1} \alpha^{2k+1}C_k + \alpha^{2l} \sum\limits_{k=l}^\infty \beta^{2k+1-2l}C_k\\ =
    \sum\limits_{k=0}^\infty \alpha^{2k+1}C_k + \alpha^{2l} \sum\limits_{k=l}^\infty \left(\beta^{2k+1-2l} - \alpha^{2k+1-2l} \right)C_k.
\end{multline}
Using \eqref{eq6.0.1} and \eqref{add_1_3} we get:
\begin{equation}
\label{add_1_4}
    u(1) \geqslant \frac{2\alpha}{1+\sqrt{1-4\alpha^2}} + \alpha^{2l}  \sum\limits_{k=0}^\infty \left(\beta^{2k+1} - \alpha^{2k+1} \right) C_{k+l}.
\end{equation}
Substituting \eqref{add_1_4} into \eqref{eigen_func_d_1} completes the proof of the lemma. \qed

\subsection*{Acknowledgments} 

S.\,A.~Molchanov came up with the idea of
the work. This study has been carried out by V.\,A.~Kutsenko and E.\,B.~Yarovaya at Steklov
Mathematical Institute of Russian Academy of Sciences. Their work was supported by the
Russian Science Foundation, project no. 23-11-00375. 

\smallskip

The authors would like to thank Prof. N.~Smorodina,  E.~Filichkina and O.~Ivlev  for helpful discussions.

\vspace{6pt} 



\begin{thebibliography}{10}

\bibitem{aizenman2015random}
M.~Aizenman and S.~Warzel, \emph{Random operators}, vol. 168, American
  Mathematical Soc., 2015.

\bibitem{akh_glazman}
N.~I. Akhiezer and I.~M. Glazman, \emph{Theory of linear operators in hilbert
  space}, Courier Corporation, 2013.

\bibitem{albeverio2000}
S.~A. Albeverio, L.~Bogachev, S.~Molchanov, and E.~Yarovaya, \emph{Annealed
  moment lyapunov exponents for a branching random walk in a homogeneous random
  branching environment}, Universit{\"a}t Bonn. SFB 256. Nichtlineare Partielle
  Differentialgleichungen, 2000.

\bibitem{anderson}
P.~W. Anderson, \emph{Absence of diffusion in certain random lattices},
  Physical review \textbf{109} (1958), no.~5, 1492.

\bibitem{catalan}
K.~N. Boyadzhiev, \emph{Series with central binomial coefficients, catalan
  numbers, and harmonic numbers}, J. Integer Seq \textbf{15} (2012), no.~1.

\bibitem{carmona2012spectral}
R.~Carmona and J.~Lacroix, \emph{Spectral theory of random schr{\"o}dinger
  operators}, Springer Science {\&} Business Media, 2012.

\bibitem{chernousova22}
E.~Chernousova, O.~Hryniv, and S.~Molchanov, \emph{Branching random walk in a
  random time-independent environment}, Mathematical Population Studies
  \textbf{30} (2023), no.~2, 73--94.

\bibitem{GM98}
J.~G{\"a}rtner and S.~Molchanov, \emph{Parabolic problems for the anderson
  model: Ii. second-order asymptotics and structure of high peaks}, Probability
  theory and related fields \textbf{111} (1998), 17--55.

\bibitem{GM90}
J.~G{\"a}rtner and S.~A. Molchanov, \emph{Parabolic problems for the anderson
  model: I. intermittency and related topics}, Communications in mathematical
  physics \textbf{132} (1990), 613--655.

\bibitem{gun15}
O.~G{\"u}n, W.~K{\"o}nig, and O.~Sekulovi{\'c}, \emph{Moment asymptotics for
  multitype branching random walks in random environment}, Journal of
  Theoretical Probability \textbf{28} (2015), no.~4, 1726--1742.

\bibitem{kato}
T.~Kato, \emph{Perturbation theory for linear operators}, vol. 132, Springer
  Science {\&} Business Media, 2013.

\bibitem{konig16}
W.~K{\"o}nig, \emph{The parabolic anderson model: Random walk in random
  potential}, Birkh{\"a}user, 2016.

\bibitem{gun13}
W.~K{\"o}nig, O.~G{\"u}n, and O.~Sekulovi{\'c}, \emph{Moment asymptotics for
  branching random walks in random environment},  (2013).

\bibitem{kunz1980spectre}
H.~Kunz and B.~Souillard, \emph{Sur le spectre des op{\'e}rateurs aux
  diff{\'e}rences finies al{\'e}atoires}, Communications in Mathematical
  Physics \textbf{78} (1980), 201--246.

\bibitem{kutsenko23}
V.~Kutsenko, D.~Sokoloff, and E.~Yarovaya, \emph{Instabilities in random media
  and peaking regimes}, Journal of Experimental and Theoretical Physics
  \textbf{136} (2023), no.~4, 498--508.

\bibitem{umn}
V.~A. Kutsenko, S.~A. Molchanov, and E.~B. Yarovaya, \emph{Supercriticality
  conditions for branching walks in a random killing medium with a single
  reproduction center}, Uspekhi matematicheskikh nauk \textbf{78} (2023), no.~5
  (473), 181--182.
\newblock In Russian.

\bibitem{M94}
S.~Molchanov, \emph{Lectures on random media}, Lectures on probability theory,
  1994.

\bibitem{shukurov84}
A.~M. Shukurov, D.~D. Sokolov, and A.~A. Ruzmaikin, \emph{Explosive growth of
  the magnetic energy in a turbulent medium.}, MHD. \textbf{19} (1984), no.~3,
  274--279.

\bibitem{yarovaya2007}
E.~Yarovaya, \emph{Branching random walks in an inhomogeneous medium}, MSU,
  Faculty of Mechanics and Mathematics \textbf{104} (2007).
\newblock In Russian.

\bibitem{yarovaya2012}
E.~Yarovaya, \emph{Symmetric branching walks in homogeneous and non homogeneous
  random environments}, Communications in Statistics-Simulation and Computation
  \textbf{41} (2012), no.~7, 1232--1249.

\bibitem{zeldovich85}
Y.~B. Zeldovich, S.~A. Molchanov, A.~A. Ruzmaikin, and D.~D. Sokolov,
  \emph{Intermittency of passive fields in random media}, Journal of
  Experimental and Theoretical Physics \textbf{89} (1985), no.~6, 2061.
\newblock In Russian.

\bibitem{zeldovich87}
Y.~B. Zeldovich, S.~A. Molchanov, A.~A. Ruzmaikin, and D.~D. Sokolov,
  \emph{Intermittency in a random environment}, Uspekhi fisicheskikh nauk
  \textbf{152} (1987), no.~5, 3--32.
\newblock In Russian.

\end{thebibliography}

\end{document}